\documentclass[reqno, 10pt]{amsart}  
\usepackage{amsmath}
\usepackage{amssymb}
\usepackage{hyperref}
\hypersetup{colorlinks=true}

\newtheorem{prop}{Proposition}[section]
\newtheorem{lemma}[prop]{Lemma}
\newtheorem{thm}[prop]{Theorem}
\newtheorem{cor}[prop]{Corollary}

\newcommand{\shrink}[1]{ {\scriptstyle {\textstyle {#1} } } }
\newcommand{\smfrac}[2]{ \shrink{ \frac{#1}{#2} } }

\newcommand{\lin}{\langle}
\newcommand{\rin}{\rangle}
\newcommand{\Sym}{\mathbb{S}^{n \times n}}
\newcommand{\Symp}{\Sym_{\plus}}

\newcommand{\plus}{{\scriptscriptstyle +}}

\newcommand{\tr}{\mathrm{ {\bf  tr}}} 
\newcommand{\sdp}{\mathrm{SDP}}

\renewcommand{\int}{\mathrm{int}}

\newcommand{\lambdamin}{\lambda_{\mathrm{min}}} 
\newcommand{\grad}{\nabla}

\begin{document}

\newpage 
$ \textrm{~} $ \quad

\title[Efficient First-Order Methods for LP and SDP]{Efficient First-Order Methods \\ for  Linear Programming \\ and Semidefinite  Programming}

\author[J. Renegar]{James Renegar  }
\address{School of Operations Research and Information Engineering,
 Cornell University, Ithaca, NY, U.S.}
 
 \date{September 19, 2014}

\thanks{Thanks to Yurii Nesterov for interesting conversations, and for encouragement, during a recent visit to Cornell  -- and special thanks for his creative research in the papers without which this one would not exist.}

\maketitle

\vspace{-8mm}

\section{{\bf  Introduction}}  \label{sect.intro}
The study of first-order methods has largely dominated research in continuous optimization for the last decade, yet still the range of problems for which ``optimal'' first-order methods have been developed is surprisingly limited, even though much has been achieved in some areas with high profile, such as compressed sensing.  Even if one restricts attention to, say, linear programming, the problems proven to be solvable by first-order methods in $ O(1/ \epsilon ) $ iterations all possess noticeably strong structure. \vspace{1mm}

We present a simple transformation of any linear program or semidefinite program  into an equivalent convex optimization problem whose only constraints are linear equations.  The objective function is defined on the whole space, making virtually all subgradient methods be immediately applicable. We observe, moreover, that the objective function is naturally ``smoothed,'' thereby allowing most first-order methods to be applied.  \vspace{1mm}

We develop complexity bounds in the unsmoothed case for a particular subgradient method, and in the smoothed case for Nesterov's original ``optimal'' first-order method for smooth functions.  We achieve the desired bounds on the number of iterations,  $ O(1/ \epsilon^2) $ and $ O(1/ \epsilon) $, respectively.    However, contrary to most of the literature on first-order methods, we measure error relatively, not absolutely.  On the other hand, also unlike most of the literature, we require only the level sets to be bounded, not the entire feasible region to be bounded.  \vspace{1mm}

Perhaps most surprising is that the transformation from a linear program or a semidefinite program is simple and so is the basic theory, and yet the approach has been overlooked until now, a blind spot. Once the transformation is realized, the remaining effort in establishing complexity bounds is mainly straightforward, by making use of various works of Nesterov. \vspace{1mm}

The following section presents the transformation and basic theory.  At the end of the section we observe that the transformation and theory extend far beyond semidefinite programming with proofs virtually identical to the ones given.  Thereafter we turn to algorithms, first for the unsmoothed case. This is where we actually rely on structure possessed by linear programs and semidefinite programs but not by conic optimization problems in general.  \vspace{1mm}

A forthcoming paper \cite{renegar2014hyperbolic}  generalizes the results to all of hyperbolic programming. That paper depends on this one.

\section{{\bf  Basic Theory}}  \label{sect.thy}

As a linear programming problem 
\[ 
 \begin{array}{rl} \min & c^T x \\
   \textrm{s.t.} & Ax = b \\
& x \geq 0 \end{array}  \]     
 is a special case -- duality aside -- of a semidefinite program in which all off-diagonal entries are constrained to equal $ 0 $, in developing the theory  we focus on semidefinite programming, as there is no point in doing proofs twice, once for linear programming and again for semidefinite programming. 
  After proving the first theorem, we digress to make certain the reader is clear on how to determine the implications of the paper for the special case of linear programming.  (We sometimes digress to consider the special case of linear programming in later sections as well.)  \vspace{1mm}

For $ C, A_1, \ldots, A_m \in \Sym $ ($ n \times n $ symmetric matrices), and $ b \in \mathbb{R}^{m} $, consider the semidefinite program   
\[ \left. \begin{array}{rl}
\inf & \lin C, X \rin  \\
\textrm{s.t.} & {\mathcal A}  (X) = b \\
  & X \succeq 0 \end{array} \right\} \sdp \]
where $ \langle \; , \; \rangle $ is the trace inner product, where  $ {\mathcal A} (X) :=  ( \lin A_1, X \rin , \ldots, \lin A_m, X \rin ) $, and where $ X \succeq 0 $ is shorthand for $ X \in \Symp $ (cone of positive semidefinite matrices). Let $ \mathrm{opt\_val} $ be the optimal value of SDP.  \vspace{1mm}

Assume $ C $ is not orthogonal to the nullspace of $ {\mathcal A} $, as otherwise all feasible points are optimal. \vspace{1mm}

Assume a strictly feasible matrix $ E $  is known.  Until section~\ref{sect.scaling}, assume $ E = I $, the identity matrix.  Assuming the identity is feasible makes the ideas and analysis particularly transparent.  In section~\ref{sect.scaling}, it is shown that the results for $ E = I $ are readily converted to results when the known feasible matrix $ E $ is a positive-definite matrix other than the identity.  Until section~\ref{sect.scaling}, however, the assumption $ E = I $ stands, but is not made explicit in the formal statement of results.    \vspace{1mm}

For symmetric matrices $ X $, let $ \lambdamin(X) $ denote the minimum eigenvalue of $ X $.  It is well known that $ X \mapsto \lambda_{\min}(X) $ is a concave function.   
\vspace{1mm}

\begin{lemma}  \label{lem.ba}
Assume $ \mathrm{SDP} $ has bounded optimal value.  If $ X \in \Sym $ satisfies $ {\mathcal A}(X) = b $ and  $ \lin C, X \rin  < \lin C, I \rin  $, then $ \lambdamin(X) < 1 \; . $ 
\end{lemma}
\noindent {\bf Proof:}  If $ \lambda_{\min}(X) \geq 1 $, then $ I + t( X - I) $ is feasible for all $ t \geq 0 $.  As the function $ t \mapsto \lin C, I + t(X- I) \rin $ is strictly decreasing (because $ \lin C, X \rin < \lin C, I \rin $), this implies SDP has unbounded optimal value, contrary to assumption. \hfill $ \Box $
 \vspace{2mm}

For all $ X \in \Sym $ for which $ \lambda_{\min} < 1 $, let $ Z(X) $ denote the matrix where the line from $ I $ in direction $ X - I $ intersects the boundary of $ \Symp $, that is,
\[  Z(X) :=  I + \smfrac{1}{1 - \lambda_{\min}(X)} (X - I)  \; . \]
We refer to $ Z(X) $ as ``the projection (from $ I $) of $ X $ to the boundary of the semidefinite cone.'' \vspace{1mm}

The following result shows that SDP is equivalent to a particular eigenvalue optimization problem for which the only constraints are linear equations.  Although the proof is straightforward, the centrality of the result to the development makes the result be a theorem. \vspace{1mm}

\begin{thm}  \label{thm.bb}
Let $ \mathrm{val} $ be any value satisfying \, $   \mathrm{val} < \lin C, I \rin  \; . $  If $ X^* $ solves
\begin{equation}  \label{eqn.ba}
  \begin{array}{rl}
   \max & \lambdamin(X) \\
 \mathrm{s.t.} &   {\mathcal A}(X) = b \\
    & \lin C, X \rin  = \mathrm{val} \; , \end{array} 
    \end{equation} 
then $ Z( X^* ) $ is optimal for $ \mathrm{SDP} $. Conversely, if $ Z^*  $ is optimal for $ \mathrm{SDP} $, then $ X^* :=   I + \frac{\lin C, I \rin - \mathrm{val} }{\lin C, I \rin - \mathrm{opt\_val}} ( Z^* - I) $ is optimal for (\ref{eqn.ba}), and $ Z^*  = Z( X^* ) $.  
\end{thm}
\noindent {\bf Proof:} Fix a value satisfying $ \mathrm{val} < \lin C, I \rin  $, and consider the affine space that forms the feasible region for (\ref{eqn.ba}):
\begin{equation}  \label{eqn.bb}
  \{ X \in \Sym: {\mathcal A}(X) = b \textrm{ and } \lin C, X \rin = \mathrm{val} \} \; . 
  \end{equation} 
Since $ \mathrm{val} < \lin C, I \rin $, it is easily proven from the convexity of $ \Symp $ that $ X \mapsto Z(X) $ gives a one-to-one map from the set (\ref{eqn.bb})   onto 
\begin{equation}  \label{eqn.bc}
  \{ Z \in \partial  \,  \Symp : {\mathcal A}(Z) = b \textrm{ and } \lin C, Z \rin < \lin C, I \rin \} \; , 
  \end{equation}
where $ \partial \,  \Symp $ denotes the boundary of $ \Symp $. \vspace{1mm}

For $ X $ in the set (\ref{eqn.bb}), the objective value of $ Z(X) $ is 
\begin{align}
  \lin C, Z(X) \rin & = \lin C, I + \smfrac{1}{1 - \lambda_{ \min}(X)} (X - I) \rin \nonumber \\
                    & = \lin C, I \rin + \smfrac{1}{1 - \lambda_{ \min }(X)} ( \mathrm{val} - \lin C, I \rin ) \; , \label{eqn.bd}
                    \end{align}
 a strictly-decreasing function of $ \lambda_{\min}(X) $.                    
 Since the map $ X \mapsto Z(X) $ is a bijection between the sets (\ref{eqn.bb}) and (\ref{eqn.bc}), solving SDP is thus equivalent to solving (\ref{eqn.ba}). \hfill $ \Box $
 \vspace{2mm}

SDP has been transformed into an equivalent linearly-constrained maximization problem with concave -- albeit nonsmooth -- objective function.  Virtually any subgradient method can be applied to this problem, the main cost per iteration being in computing a subgradient and projecting it onto the subspace $ \{ V: {\mathcal A}(V) = 0 \textrm{ and } \lin C, V \rin = 0 \} $.    In section \ref{sect.nesterov}, it is observed that the objective function has a natural smoothing, allowing almost all first-order methods to be applied, not just subgradient methods.   
\vspace{1mm}

We digress to interpret the implications of the development thus far for the linear programming problem
\begin{equation}  \label{eqn.bz}
  \left. \begin{array}{rl}
 \min & c^T x \\
\textrm{s.t.} & Ax = b \\
 & x \geq 0 \; . \end{array}   \right\} \, \mathrm{LP}  
 \end{equation} 
LP can easily be expressed as a semidefinite program in variable $ X \in \Sym $ constrained to have all off-diagonal entries equal to zero, where the diagonal entries correspond to the original variables $ x_1, \ldots, x_n $.  \vspace{1mm}

In particular, the standing assumption that $ I $ is feasible for SDP becomes, in the special case of LP, a standing assumption that $ \mathbf{1}  $ (the vector of all ones) is feasible.  
The eigenvalues of $ X $ become  the coordinates $ x_1, \ldots, x_n $.  The map $ X \mapsto \lambda_{\min}(X) $ becomes $ x \mapsto \min_j x_j $.  Lemma~\ref{lem.ba}  becomes the statement that if $ x $ satisfies $ Ax = b $ and $ c^T x < c^T \mathbf{1} $, then $ \min_j x_j < 1 $.  \vspace{1mm}

Finally, Theorem~\ref{thm.bb} becomes the result that for any value satisfying $ \mathrm{val} < c^T \mathbf{1}  $, LP is equivalent to
\begin{equation} \label{eqn.bda}
  \begin{array}{rl}
    \max_x & \min_j x_j \\
     \textrm{s.t.} & Ax = b \\
           & c^T x = \mathrm{val} \; , \end{array} \end{equation} 
in that, for example, $ x^* $ is optimal for (\ref{eqn.bda}) if and only if the projection $ z(x^*) = \mathbf{1}  + \smfrac{1}{1 - \min_j x^*_j} (x^* - \mathbf{1} ) $ is optimal for LP. \vspace{1mm}

In this straightforward manner, the reader can realize the implications for LP of all results in the paper. \vspace{1mm}

Before leaving the simple setting of linear programming, we make observations pertinent to applying subgradient methods to solving (\ref{eqn.bda}), the problem equivalent to LP.   \vspace{1mm}

The subgradients of $ x \mapsto \min_j x_j $ at $ x $ are  the convex combinations of the standard basis vectors $ e(k) $ for which $ x_k = \min_j x_j $. Consequently, the projected subgradients at $ x $ are the convex combinations of the vectors $ \bar{P}_k $ for which $ x_k = \min_j x_j $, where $ \bar{P}_k $ is the $ k^{th} $ column of the matrix projecting $ \mathbb{R}^n $ onto  the nullspace of $ \bar{A} = \left[ \begin{smallmatrix}  A \\ c^T \end{smallmatrix} \right] $, that is  
\begin{equation}  \label{eqn.bdb}
   \bar{P}  := I - \bar{A}^T (\bar{A} \,   \bar{A}^T)^{-1} \bar{A}  \; . 
   \end{equation} 
In particular, if for a subgradient method the current iterate is $ x $, then the chosen projected subgradient can simply be any of the columns $ \bar{P}_k $ for which $ x_k = \min_j x_j $. Choosing the projected subgradient in this way gives the subgradient method a combinatorial feel.  If, additionally, the subgradient method does exact line searches, then the algorithm possesses distinct combinatorial structure.  (In this regard it should be noted that the work required for an exact line search is only $ O( n \log n) $, dominated by the cost of sorting.) \vspace{1mm}

If $ m \ll n $, then $ \bar{P} $ is not computed in its entirety, but instead the matrix $ \bar{M}  = ( \bar{A} \bar{A}^T)^{-1} $ if formed as a preprocessing step, at cost $ O(m^2 n) $. Then, for any iterate $ x $ and an index $ k $ satisfying $ x_k = \min_j x_j $, the projected  subgradient $ \bar{P}_k $ is computed according to
 \[  u = \bar{M} \, \bar{A}_k \quad \rightarrow \quad  v = \bar{A}^T u \quad \rightarrow \quad \bar{P}_k = e(k) - v \; , \]
for a cost of $ \, O(m^2 \, + \, \# \mathrm{non\_zero\_entries\_in\_} A   \,  +  \, n  \log n) \, $ per iteration, where $ O( n  \log n) $ is the cost of finding a smallest coordinate of $ x $.  \vspace{1mm}

Now we return to the theory, expressed for SDP, but interpretable for LP in the straightforward manner explained above. \vspace{1mm}

Assume, henceforth, that SDP has at least one optimal solution. Thus,  the equivalent problem (\ref{eqn.ba}) has at least one optimal solution.     Let $ X^*_{\mathrm{val}} $ denote any of the optimal solutions for the equivalent problem.
\vspace{1mm}

\begin{lemma}  \label{lem.bc}
 \[   \lambda_{\min}(X^*_{\mathrm{val} }) = \frac{\mathrm{val} - \mathrm{opt\_val}}{\lin C, I \rin - \mathrm{opt\_val}} \] 
  \end{lemma}
\noindent {\bf Proof:}   
By Theorem~\ref{thm.bb}, $ Z(X^*_{\mathrm{val}}) $ is optimal for SDP -- in particular, $ \lin C, Z(X^*_{\mathrm{val}}) \rin = \mathrm{opt\_val} $.   Thus, according to (\ref{eqn.bd}),
\[  \mathrm{opt\_val} = \lin C, I \rin + \smfrac{1}{1 - \lambda_{\min}(X^*_{\mathrm{val} })} \, ( \mathrm{val} - \lin C, I \rin) \; . \]
Rearrangement completes the proof.      \hfill $ \Box $
 \vspace{3mm}

We focus on the goal of computing a matrix $ Z $ that is feasible for SDP and has objective value which is significantly better than the objective value for $ I $, in the sense that 
\begin{equation}  \label{eqn.be}
     \frac{\lin C, Z \rin - \mathrm{opt\_val}}{\lin C,I \rin - \mathrm{opt\_val}} \leq \epsilon \; , 
     \end{equation} 
where $ \epsilon > 0 $ is user-chosen.  Thus, for the problem of main interest, SDP (or the special case, LP), the focus is on relative improvement in the objective value. 
\vspace{1mm}

 The following proposition provides a useful characterization of the accuracy needed in approximately solving the SDP equivalent problem (\ref{eqn.ba}) so as to ensure that for the computed matrix $ X $, the projection $ Z = Z(X) $ satisfies (\ref{eqn.be}). \vspace{1mm}

\begin{prop}  \label{prop.bd}
Let $ 0 \leq  \epsilon < 1  $, and let $ \mathrm{val} $ be a value satisfying $ \mathrm{val} < \lin C, I \rin \; . $   \vspace{1mm}

If $ X $ is feasible for the $ \mathrm{SDP} $  equivalent problem (\ref{eqn.ba}), then 
\begin{gather} 
    \frac{\lin C, Z(X) \rin - \mathrm{opt\_val}}{\lin C,I \rin - \mathrm{opt\_val}} \, \leq \,  \epsilon    \label{eqn.bf} \\   \textrm{if and only if}  \nonumber  \\
    \lambdamin( X^*_{ \mathrm{val}}  ) - \lambdamin( X) \, \leq \, \frac{\epsilon }{1 - \epsilon } \, \, \frac{\lin C, I \rin - \mathrm{val} }{\lin C, I \rin - \mathrm{opt\_val} } \; .  \label{eqn.bg}
\end{gather} 
\end{prop}
\noindent {\bf Proof:}  Assume $ X $ is feasible for the equivalent problem (\ref{eqn.ba}).  For $ Y = X, X^*_{\mathrm{val}} \; , $  we have the equality (\ref{eqn.bd}), that is,  
\[     \lin C, Z(Y) \rin = \lin C, I \rin + \smfrac{1}{1 - \lambda_{\min}(Y)} ( \mathrm{val} - \lin C, I \rin ) \; . 
\]
Thus,
\begin{align*}
    \frac{\lin C, Z(X) \rin - \mathrm{opt\_val}}{\lin C,I \rin - \mathrm{opt\_val}} & =   \frac{\lin C, Z(X) \rin - \lin C, Z(X^*_{\mathrm{val}} ) \rin }{\lin C,I \rin - \lin C, Z(X^*_{\mathrm{val}} ) \rin } \\
    & = \frac{ \smfrac{1}{1 - \lambdamin(X)} -  \smfrac{1}{1 - \lambdamin(X^*_{\mathrm{val}} )} }{ - \smfrac{1}{1 - \lambdamin(X^*_{\mathrm{val}} )}} \\  
    & = \frac{ \lambdamin(X^*_{\mathrm{val}} ) - \lambdamin(X) }{ 1 - \lambdamin(X)} \; . 
    \end{align*}
Hence,
\begin{gather*}     
 \frac{\lin C, Z(X) \rin - \mathrm{opt\_val}}{\lin C,I \rin - \mathrm{opt\_val}} \leq \epsilon \\
\Leftrightarrow \\ 
   \lambdamin(X^*_{\mathrm{val}} ) - \lambdamin(X) \leq \epsilon \, ( 1 - \lambdamin(X)) \\
\Leftrightarrow \\
(1 - \epsilon ) ( \lambdamin(X^*_{\mathrm{val}} ) - \lambdamin(X) ) \leq \epsilon ( 1 - \lambdamin(X^*_{\mathrm{val}} )) \\
\Leftrightarrow \\ \lambdamin(X^*_{\mathrm{val}} ) - \lambdamin(X) \leq \smfrac{\epsilon }{1 - \epsilon } ( 1 - \lambdamin(X^*_{\mathrm{val}} )) \; .  
   \end{gather*}
Using Lemma~\ref{lem.bc}   to substitute for the rightmost occurrence of $ \lambda_{\min}(X^*_{\mathrm{val}} ) $ completes the proof. \hfill $ \Box $
 \vspace{3mm}
 
It might seem that to make use in complexity analysis of the equivalence of (\ref{eqn.bg})  with (\ref{eqn.bf}), it would be necessary to assume as input to algorithms a lower bound on $ \mathrm{opt\_val} $.  Such is not the case, as is shown in the following sections. \vspace{1mm}

In concluding the section, we observe that the basic theory holds far more generally.  In particular, let $ K $ be a closed, pointed, convex cone in $ \mathbb{R}^n $, and assume $ e $ lies in the interior of $ K $.  For $ x \in \mathbb{R}^n $, define
\begin{equation}  \label{eqn.bga}
    \lambda_{\mathrm{\min},e}(x) = \inf \{ \lambda \in \mathbb{R}: x - \lambda e \notin K \} \; . 
    \end{equation} 
It is easy to show $ x \mapsto \lambda_{\min,e}(x) $ is a closed, concave function with finite value for all $ x \in  \mathbb{R}^n $.  
\vspace{1mm}

Assume additionally that $ e $ is feasible for the conic optimization problem 
\begin{equation}  \label{eqn.bh}
   \begin{array}{rl}
    \min & \lin c, x \rin \\
   \textrm{s.t.} & Ax = b \\
    & x \in K  \end{array} \end{equation} 
    (where $ \langle \; , \; \rangle $ is any fixed inner product).
The same proof as for Theorem~\ref{thm.bb} then shows that for any value satisfying $ \mathrm{val} < \lin c, e \rin $,     
 (\ref{eqn.bh}) is equivalent to the linearly-constrained optimization problem
\begin{equation}  \label{eqn.bi}
     \begin{array}{rl} \max & \lambda_{ \min,e }(x) \\
\textrm{s.t.} & Ax = b \\
&   \lin c, x \rin = \mathrm{val} \; , \end{array} 
\end{equation} 
in the sense that if $ x $ is optimal for (\ref{eqn.bi}), then $ z(x) := e + \smfrac{1}{1 - \lambda_{\min,e} (x) } (x - e ) $ is optimal for (\ref{eqn.bh}), and conversely, if $ z^* $ is optimal for (\ref{eqn.bh}), then $ x^* = e + \smfrac{\lin c, e \rin - \mathrm{val} }{\lin c, e \rin - \mathrm{opt\_val} } (z^* - e) $ is optimal for (\ref{eqn.bi}); moreover, $ z^* = z(x^*) $.    \vspace{1mm}

Likewise, Proposition~\ref{prop.bd} carries over with the same proof.    \vspace{1mm}

Furthermore, analogous results are readily developed for a variety of different forms of conic optimization problems.  Consider, for example, a problem
\begin{equation}  \label{eqn.bj}
   \begin{array}{rl}
      \min & \lin c, x \rin \\
        \textrm{s.t.} & Ax - b \in K \; . \end{array} \end{equation} 
Now assume known a feasible point $ e'  $ in the interior of the feasible region.  Then, for any value satisfying $ \mathrm{val} < \lin c, e'  \rin $,  the conic optimization problem (\ref{eqn.bj}) is equivalent to a problem for which there is only one linear constraint:
\begin{equation}  \label{eqn.bk}
\begin{array}{rl}
   \max & \lambda_{\min, e} (Ax - b) \\
   \textrm{s.t.} & \lin c, x \rin = \mathrm{val} \; , \end{array} 
   \end{equation}
where $ e := A e'  - b $, and where $ \lambda_{\min ,e} $ is as  defined in (\ref{eqn.bga}). The problems are equivalent in that if $ x^* $ is optimal for (\ref{eqn.bk}), then the projection $ z(x^*) $ of $ x^* $ from $ e' $ to the boundary of the feasible region is optimal for (\ref{eqn.bj})  -- that is, 
\[     z(x^*) := e' + \smfrac{1}{1 - \lambda_{\min, e} (A x^* - b) } \,  (x^* - e' )   \]
is optimal for (\ref{eqn.bj}) -- and, conversely, if $ z^* $ is optimal for (\ref{eqn.bj}), then $ x^* = e'  + \smfrac{\lin c, e' \rin - \mathrm{val} }{\lin c, e' \rin - \mathrm{opt\_val}} (z^* - e') $ is optimal for  (\ref{eqn.bk}); moreover, $ z^* = z(x^*) $. \vspace{2mm}

We focus on the concrete setting of semidefinite programming (and linear programming) because the algebraic structure thereby provided is sufficient for designing provably-efficient first-order methods, in both smoothed and unsmoothed settings. We now begin validating the claim. \vspace{1mm}

\section{ {\bf Corollaries for a Subgradient Method}} \label{sect.subgrad}

Continue to assume SDP has an optimal solution and $ I $ is feasible. 
\vspace{1mm}

  Given $ \epsilon > 0 $ and a value satisfying $ \mathrm{val} < \lin C, I \rin $, we wish to approximately solve the SDP equivalent problem
\begin{equation} \label{eqn.ca}
 \begin{array}{rl}
 \max & \lambda_{ \min}(X) \\
\textrm{s.t.} & {\mathcal A}(X) = b \\
      & \lin C, X \rin = \mathrm{val} \; , \end{array} 
      \end{equation}
where by ``approximately solve'' we mean that feasible $ X $ is computed for which 
\[ 
\lambda_{\min}(X^*_{\mathrm{val}}) - \lambda_{ \min}(X) \leq \epsilon' \;  
\] 
with $ \epsilon' $ satisfying
\[ 
 \epsilon' \, \leq \, \frac{\epsilon}{1 - \epsilon } \, \, \frac{\lin C, I \rin - \mathrm{val}   }{\lin C, I \rin - \mathrm{opt\_val} }  \; . 
\] 
Indeed, according to Proposition~\ref{prop.bd}, the projection $ Z = Z(X) $ will then satisfy
\begin{equation}  \label{eqn.cb}
\frac{\lin C, Z \rin - \mathrm{opt\_val}}{\lin C,I \rin - \mathrm{opt\_val}} \leq \epsilon  \; . 
\end{equation}
\vspace{1mm}

We begin by recalling a well-known complexity result for a subgradient method, interpreted for when the method is applied to solving the SDP equivalent problem (\ref{eqn.ca}).  From this is deduced a bound on the number of iterations sufficient to obtain $ X $ whose projection $ Z = Z(X) $ satisfies (\ref{eqn.cb}). We observe, however, that in a certain respect, the result is  disappointing.  In the next section, the framework  is embellished by applying the subgradient method not to  (\ref{eqn.ca}) for only one value $ \mathrm{val} $, but to (\ref{eqn.ca}) for  a small and carefully chosen sequence of values, $ \mathrm{val} = \mathrm{val}_{\ell } $.   The embellishment results in a computational scheme which possesses the desired improvement. \vspace{1mm}

For specifying a subgradient method and stating a bound on its complexity, we follow Nesterov's book \cite{nesterov2004introductory}:
\begin{itemize}

\item {\bf  Subgradient Method}
\begin{enumerate}

\addtocounter{enumi}{-1}

\item Inputs:
\begin{itemize}

\item Number of iterations: $ N $ 
\item Initial iterate: $ X_0 $ satisfying $ {\mathcal A}(X_0) = b $ and $ \lin C, X_0 \rin < \lin C, I \rin. $  \\
$ \textrm{~} $ \qquad \qquad  \qquad  Let $ \mathrm{val} := \lin C, X_0 \rin $. 
\item Distance upper bound: $ R $, a value for which there exists $ X^*_{\mathrm{val}} $ \\
$ \textrm{~} $ \qquad \qquad  \qquad  \qquad  \qquad  \qquad    satisfying $ \| X_0 - X^*_{\mathrm{val}}  \| \leq R \; . $   
 \item Initial ``best'' iterate: $ X = X_0 $
 \item  Initial counter value: $ k = -1 $  
 \end{itemize} 
\item Update counter: $ k+1 \rightarrow k $ 
\item Iteration:
Compute a subgradient $ \grad \lambda_{\min}(X_k) $ and orthogonally project it onto the subspace
\begin{equation}  \label{eqn.cc}
    \{ V: {\mathcal A}(V) = 0 \textrm{ and } \lin C, V \rin = 0 \} \; . \end{equation} 
Denoting the projection by $ G_k $, compute
\[     X_{k+1} := X_k + \smfrac{R}{\sqrt{N} \,  \| G_k \| } G_k \; .\]
\item If $ \lambda_{\min}(X_{k+1}) > \lambda_{\min}(X) $,  then make the replacement $ X_{k+1} \rightarrow X \; . $
\item Check for termination: If $ k = N - 1 $,  then output $ X $ and terminate. \\
$ \textrm{~} $ \qquad \qquad  \qquad  \qquad  \qquad    Else, go to Step 1. 
\end{enumerate}
\end{itemize}
\vspace{1mm}

\begin{thm}  \label{thm.ca}
For $ \mathrm{Subgradient \, \,  Method} $, the output $ X  $ satisfies
\[ 
          \lambda_{\min}(X^*_{\mathrm{val}}) - \lambda_{\min}( X) \, \leq \, R/ \sqrt{N}   \; , 
          \]
where $ \mathrm{val} := \lin C, X_0 \rin $ and $ X_0 $ is the input matrix.       
\end{thm}
\noindent {\bf Proof:}  The function $ X \mapsto \lambda_{\min}(X) $ is  Lipschitz continuous with constant 1:
\[  | \lambda_{\min}(X) - \lambda_{\min}(Y) | \leq \| X - Y \| \quad = \big(  \sum_j \lambda_j ( X - Y)^2 \big)^{1/2}  \; . \]
The result is thus a simple corollary of Theorem 3.2.2 in Nesterov's book, by choosing the parameter values there to be $ h_k = R/\sqrt{N} $ for $ k = 0, \ldots, N-1 $.    \vspace{3mm} \hfill $ \Box $

 We briefly digress to the special case of linear programming.  \vspace{1mm}

  Recall for LP  -- the linear program (\ref{eqn.bz}) --  the projected subgradients at $ x $ are the convex combinations of the columns $ \bar{P}_k $ for which $ x_k = \min_j x_j $, where $ \bar{P} $ is the matrix projecting $ \mathbb{R}^n $ orthogonally onto the nullspace of $ \bar{A} = \left[ \begin{smallmatrix} A \\ c^T \end{smallmatrix} \right] $, that is,
\[  \bar{P} = I - \bar{A}^T (\bar{A} \, \bar{A}^T)^{-1} \bar{A} \; . \]
         In particular, if $ x $ is the current iterate for Subgradient Method, a projected subgradient can be selected simply by computing any column $ \bar{P}_k $   for which $ x_k = \min_j x_j $.  Subgradient Method then moves from $ x $ to $ x + \frac{R}{\sqrt{N} \, \| \bar{P}_k \|} \bar{P}_k $.  The geometry is interesting in that each step is being chosen from among only the vectors $ \frac{R}{\sqrt{N} \, \| \bar{P}_j \|} \bar{P}_j $ for $ j = 1, \ldots ,n $. \vspace{1mm}

  The geometry is made even more interesting by Theorem~\ref{thm.ca} asserting that even for the choice of steps coming from this limited set of vectors, still it holds that the final output $ x $ satisfies
\[  \min_j x^*_{\mathrm{val},j} \, - \, \min_j x_j \leq R/ \sqrt{N} \; , \]
where $ x^*_{\mathrm{val},j} $ denotes the $ j^{th} $ coordinate of the optimal solution $ x^* $ for the LP equivalent problem
\[ \begin{array}{rl}
          \max_x & \min_j x_j \\
          \textrm{s.t.} & Ax = b \\
                 & c^T x = \mathrm{val} \; . 
                 \end{array} \]     
\vspace{2mm}

Now we return to the more general setting of semidefinite programming. \vspace{1mm}

Below, the input matrix $ X_0 $ to Subgradient Method is required to be feasible for SDP, mainly so that the input $ R $ can be chosen as a value with clear relevance to SDP, a value we now describe. \vspace{1mm}

 The ``level sets'' for SDP are the sets 
\[  \mathrm{Level}_{\textrm{val}} =   \{ X \succeq 0 : {\mathcal A}(X) = b \textrm{ and } \lin C, X \rin  = \textrm{val} \} \; ,  \]
where $ \textrm{val} $ is any fixed value. Of course $ \mathrm{Level}_{\mathrm{val}} = \emptyset $ if $ \mathrm{val} < \mathrm{opt\_val} $. \vspace{1mm}

     If some nonempty level set is bounded, then all level sets are bounded. On the other hand, if a level set is unbounded, then either SDP has unbounded optimal value or can be made to have unbounded value with an arbitrarily small perturbation of $ C $.  Thus, in developing numerical methods for approximating optimal solutions, it is natural to focus on the case that level sets for SDP are bounded (equivalently, the dual problem is strictly feasible). Hence, we assume the level sets are bounded.   \vspace{1mm}  
     
Let $ \mathrm{diam} $ be a known value satisfying
\[  \mathrm{diam} \geq \max \{ \|X - Y \|: X, Y \in \mathrm{Level}_{\mathrm{val}} \} \quad  \textrm{ for all $ \mathrm{val} < \lin C, I \rin $}  \; ,  \]
that is, an upper bound on the diameters of all level sets for better objective values than the value for the level set containing $ I $.   \vspace{1mm}

  Although the assumption of knowing the upper bound $ \mathrm{diam} $ is strong, it is consistent with assumptions found throughout the literature on first-order methods, such as the requirement for Subgradient Method that the input $ R $ be an upper bound on $ \| X_0 - X^*_{\mathrm{val}} \| $, where $ X_0 $ is the input matrix.  \vspace{1mm}

Moreover, even though the assumption of knowing $ \mathrm{diam} $ is strong,    still there are many interesting situations in which the assumption is valid, particularly when a problem is specifically modeled in such a way as to make the diameter of the level sets (for $ \mathrm{val} < \lin C, I \rin $) be of reasonable magnitude.  For example, when $ I $ is on (or near) the central path, the choice of upper bound $ \mathrm{diam} = n $ is valid, albeit for various carefully modeled semidefinite programs in which $ I $ is explicitly made to be near the central path, stronger upper bounds hold (e.g., $ \mathrm{diam} = O(  \sqrt{n}) $ in numerous interesting cases, some of which are displayed in the forthcoming paper \cite{renegar2014hyperbolic}).  \vspace{1mm}
  
In most of the literature on optimal first-order methods, the feasible region is required to be bounded, not just the level sets.  By focusing on relative error (\ref{eqn.cb}) rather than absolute error, we are able to require only that the level sets be bounded, not the feasible region.   \vspace{1mm}

In the following corollary regarding Subgradient Method, the choice of input $ N $ depends on the optimal value for SDP.  Naturally the reader will infer that in addition to knowing the upper bound $ \mathrm{diam} $, our algorithmic scheme will require knowing a lower bound on $ \mathrm{opt\_val} $, but this is not the case for the scheme.  The corollary is used for motivating the next step in specifying the scheme. \vspace{1mm}

\begin{cor} \label{cor.cb} Assume $ X_0 $ is feasible for $ \mathrm{SDP} $ and satisfies $ \lin C, X_0 \rin < \lin C, I \rin $.  Define $ \mathrm{val} :=  \lin C, X_0 \rin $.  Let $ 0 < \epsilon < 1 $. \vspace{1mm}

 If $ X_0 $ and $ R = \mathrm{diam} $ are inputs to $ \mathrm{Subgradient \, \,  Method} $, along with an integer $ N $ satisfying
\begin{equation} \label{eqn.ce}
     N \, \geq \, \left(    \frac{\mathrm{diam}}{ \epsilon} \, \, \, \frac{\lin C, I \rin - \mathrm{opt\_val}}{\lin C, I \rin - \mathrm{val} } \right)^2  \; , 
\end{equation} 
then for the output $ X $, the projection $ Z(X) $  satisfies
\[       \frac{\lin C, Z(X) \rin - \mathrm{opt\_val}}{\lin C, I \rin - \mathrm{opt\_val}} \, \leq \, \epsilon \; . \]
\end{cor}
\noindent {\bf Proof:}  Since $ X_0 $ is feasible, so is $ X^*_{\mathrm{val}} $ (because $ 0 \leq \lambda_{\min}(X_0) \leq \lambda_{\min}(X^*_{\mathrm{val}}) $).  Thus, $ \| X_0 - X^*_{\mathrm{val}}  \| \leq \mathrm{diam} $, making $ R = \mathrm{diam} $ a valid input to Subgradient Method. \vspace{1mm}

For inputs as specified, Theorem~\ref{thm.ca}   immediately implies the output $ X $ for Subgradient Method satisfies
\[       \lambda_{\min}( X^*_{\mathrm{val}}) - \lambda_{\min}(X) \, \leq \, \epsilon \, \cdot \,  \frac{\lin C, I \rin - \mathrm{val}}{ \lin C, I \rin - \mathrm{opt\_val}} \; . \]  
Invoking Proposition~\ref{prop.bd}  completes the proof.    \vspace{3mm}  \hfill $ \Box $

The dependence of the iteration lower bound (\ref{eqn.ce})  on $ \epsilon^2 $ is unfortunate but probably unavoidable without smoothing the objective function $ X \mapsto \lambda_{\min}(X) $, as is done in sections~\ref{sect.nesterov} and \ref{sect.smoothed}.  Likewise, a significant dependence on $ \mathrm{diam} $ -- or some other meaningful quantity capturing the distance $ \| X_0 - X^*_{\mathrm{val}}  \| $ -- probably is unavoidable.  However, the dependence on $  \frac{\lin C, I \rin - \mathrm{opt\_val}}{\lin C, I \rin - \mathrm{val} } $ is  disconcerting. \vspace{1mm}

To understand why the dependence is disconcerting, consider that the most natural choice for the input matrix is $ X_0 = I - \smfrac{1}{\lambda_{\max}( \pi(C)) } \pi(C) $, where $ \pi(C) $ is the orthogonal projection of $ C $ onto the subspace $ \{ V: {\mathcal A}(V) = 0 \} $.  This is the choice for $ X_0 $ obtained by moving from $ I $ in direction $ - \pi(C) $ until the boundary of the semidefinite cone is reached. \vspace{1mm}

   However, {\em even when $ I $ is on the central path} (in which case the direction $ -\pi(C) $ is tangent to the central path), it can happen that the value $  \frac{\lin C, I \rin - \mathrm{opt\_val}}{\lin C, I \rin - \mathrm{val} } $ is of magnitude $ \sqrt{n}  $ for $ \mathrm{val} = \lin C, X_0 $ and $ X_0 = I - \smfrac{1}{\lambda_{\max}(\pi(C))} \pi(C) $.  Thus, even for problems modeled carefully so that $ I $ is on  the central path and $ \mathrm{diam} $ is of limited size, the iteration lower bound (\ref{eqn.ce}) can grow significantly with $ n $  regardless of the value for $ \epsilon $.  This is disconcerting. \vspace{1mm}

Moreover, we want an algorithm for which $ \mathrm{opt\_val} $ does not explicitly figure into choosing the inputs.  We already assume the upper bound $ \mathrm{diam} $ is known.  We want to avoid also assuming a lower bound on $ \mathrm{opt\_val} $ is known.  \vspace{1mm}

These matters are handled in the following section.  

\section{{\bf The NonSmoothed Scheme}}  \label{sect.nonsmoothed}

The observations concluding the preceding section raise a question: 
\begin{quote}
 Is it possible to efficiently move from an initial feasible matrix $ U_0 $ satisfying $ \lin C, U_0 \rin < \lin C, I \rin $, to a feasible matrix $ Y $ for which $ \mathrm{val} = \lin C, Y \rin $ satisfies, say, $ \smfrac{\lin C, I \rin - \mathrm{opt\_val} }{\lin C, I \rin - \mathrm{val}  } \leq 3 $?
 \end{quote}
   We begin this section by providing an affirmative answer, but first let us again display the pertinent optimization problem:
\begin{equation} \label{eqn.da}
\begin{array}{rl}
 \max & \lambda_{ \min}(X) \\
\textrm{s.t.} & {\mathcal A}(X) = b \\
   & \lin C, X \rin = \mathrm{val} \; . 
   \end{array} \end{equation}
   Recall that  $ X^*_{ \mathrm{val}} $ denotes any optimal solution of (\ref{eqn.da}), an optimization problem which is equivalent to SDP (assuming $ \mathrm{val} < \lin C, I \rin $). 
   \vspace{1mm}

Consider the following computational procedure: 
\begin{itemize}

\item {\bf NonSmoothed SubScheme}
\begin{enumerate}

\addtocounter{enumi}{-1}

\item Initiation:

\begin{itemize}

\item Input:  A matrix $ U_0 $ that is feasible for SDP and satisfies $ \lin C, U_0 \rin < \lin C, I \rin $.

\item Let $ \mathrm{val}_0 = \lin C, U_0 \rin $ 
\item Let $ \ell = -1 $. 
\end{itemize}

\item Outer Iteration Counter Step: $ \ell + 1 \rightarrow \ell $

\item  Inner Iterations:

\begin{itemize}
\item Apply Subgradient Method with inputs $ X_0 = U_{\ell} $, $ R = \mathrm{diam} $ and   $ N = \lceil 9 \, \mathrm{diam}^2  \rceil $.
\item  Rename the output $ X $ as  $ V_{ \ell} $.  

\end{itemize}

\item Check for Termination: 

\begin{itemize}

\item If $ \lambda_{\min}( V_{\ell}) \leq 1/3 $, then output $ Y = U_{ \ell} $ and terminate.  

\item Else, compute the projection
\[  \qquad \qquad  \qquad  \qquad   U_{ \ell + 1} := Z(V_{\ell}) \; , \quad  \textrm{let  } \mathrm{val}_{\ell + 1} := \lin C, U_{ \ell + 1} \rin \; ,   \]
and go to Step 1.     
 
\end{itemize}

\end{enumerate}
\end{itemize}
\vspace{1mm}

\begin{prop} \label{prop.da}
$ \mathrm{NonSmoothed \, \,  SubScheme} $  outputs $ Y $ that is feasible for SDP and satisfies 
\[  \frac{\lin C, I \rin - \mathrm{opt\_val} } {\lin C, I \rin - \mathrm{val}  }\,  \leq \, 3  \; , \]
  where $ \mathrm{val} := \lin C, Y \rin $.  The total number of outer iterations does not exceed 
\[  \log_{3/2} \left( \frac{\lin C, I \rin - \mathrm{opt\_val} }{\lin C, I \rin - \mathrm{val}_0 } \right)  \; ,      \]
where $ \mathrm{val}_0 = \lin C, U_0 \rin $ and $ U_0 $ is the input matrix. 
\end{prop}
\noindent {\bf Proof:}  It is easily verified that all of the matrices $ U_{ \ell} $ and $ V_{ \ell} $ computed by NonSmooth SubScheme satisfy the SDP equations $ {\mathcal A}(X) = b $.  Moreover, $ U_{ \ell} $ is clearly feasible for SDP, lying in the boundary of the feasible region.   \vspace{3mm}

Fix $ \ell $ to be any value attained by the counter.  We now examine the effects of Steps 2 and 3. \vspace{1mm}

Corollary~\ref{cor.cb}  with $ N = \lceil 9 \, \mathrm{diam}^2   \rceil  $ shows that in Step 2, the output $ X $ from Subgradient Method satisfies 
 \[  \lambda_{\min}(X^*_{\mathrm{val}_{\ell}}) -  \lambda_{\min}(X) \leq 1/3 \; , \]
that is, 
\begin{equation}   \label{eqn.db}
    \lambda_{\min}(X^*_{\mathrm{val}_{\ell}}) -  \lambda_{\min}(V_{\ell}) \leq 1/3 \; . 
    \end{equation}
    
Observe 
\begin{align*} 
\frac{\lin C, I \rin - \mathrm{val}_{ \ell} }{\lin C, I \rin - \mathrm{opt\_val}} & = 1 -  \lambda_{\min}(X^*_{\mathrm{val}_{\ell}}) \quad \textrm{(by Lemma~\ref{lem.bc})}  \\
  & \geq \smfrac{2}{3} - \lambda_{\min}( V_{ \ell}) \quad \textrm{(by (\ref{eqn.db}))}
  \end{align*}
Hence, if the method terminates in Step 3 -- that is, if $ \lambda_{\min}(V_{\ell}) \leq 1/3 $ -- then the output matrix $ Y = U_{ \ell} $ satisfies
\begin{equation}  \label{eqn.dc}
    \frac{\lin C, I \rin - \mathrm{val} }{\lin C, I \rin - \mathrm{opt\_val}} \, \geq \, \frac{1}{3} \; , 
    \end{equation} 
where $ \mathrm{val} :=  \lin C, Y \rin  = \mathrm{val}_{ \ell} $.  We have now verified that if NonSmoothed SubScheme terminates, then  the output $ Y $ is indeed feasible for SDP and satisfies the desired inequality (\ref{eqn.dc}).  \vspace{1mm}

On the other hand, if the method does not terminate in Step 3, it computes the matrix $ U_{ \ell + 1} $ and its objective value, $ \mathrm{val}_{ \ell + 1} $.  Here, observe 
\begin{align*}
 \mathrm{val}_{ \ell + 1 }  & = \lin C, I \rin +  \smfrac{1}{1 - \lambda_{\min}(V_{\ell})} \, ( \mathrm{val}_{\ell} - \lin C, I \rin ) \\
  & \leq \lin C, I \rin -   \smfrac{3}{2}  ( \lin C, I \rin - \mathrm{val}_{\ell} )  \; , 
  \end{align*}
because $ \mathrm{val}_{\ell} < \lin C, I \rin $ and  $ \lambda_{\min}(V_{ \ell}) \geq 1/3 $ (due to no termination). 
Hence,
\[   \frac{\lin C, I \rin - \mathrm{val}_{ \ell + 1} }{\lin C, I \rin - \mathrm{opt\_val}}  \geq \frac{3}{2} \, \, \frac{\lin C, I \rin - \mathrm{val}_{ \ell} }{\lin C, I \rin - \mathrm{opt\_val}} \; . \]
Since all values $ \mathrm{val}_{\ell} $ computed by the algorithm satisfy $ \mathrm{val}_{\ell} \geq \mathrm{opt\_val} $ (as $ U_{ \ell } $ is feasible for SDP),  it immediately follows that 
\[  \log_{3/2} \left(  \frac{\lin C, I \rin - \mathrm{opt\_val}} {\lin C, I \rin - \mathrm{val}_0 }  \right) \]
is an upper bound on the number of outer iterations.  \vspace{3mm} \hfill $ \Box $

 Specifying our overall computational scheme relying on the subgradient method, and analyzing the scheme's  complexity, both are now easily accomplished:
\begin{itemize}

\item {\bf NonSmoothed Scheme}

\begin{enumerate}

\addtocounter{enumi}{-1}

\item Inputs: A value $ 0 < \epsilon < 1 $, and a matrix $ U_0 $ which both is feasible for $ \mathrm{SDP} $ and satisfies $ \lin C, U_0 \rin < \lin C, I \rin $.  

\begin{itemize}

\item For example, the matrix $ U_0 = I -  \smfrac{1}{\lambda_{\max}(\pi(C))} \, \pi(C) $.

\end{itemize}
 
\item Apply NonSmoothed SubScheme with input $ U_0 $.  Let $ Y $ denote the output.

\item Apply Subgradient Method with inputs $ X_0 = Y $, $ R = \mathrm{diam} $ and 
       \[  N = \left\lceil \left( 3 \, \mathrm{diam}/\epsilon  \right)^2   \right\rceil \; . \]
   Let $ X $ denote the output.  

\item Compute and output the projection $ Z = Z(X) $, then terminate.   

\end{enumerate}

\end{itemize}
\vspace{2mm}

In stating the following theorem, we make explicit that $ I $ is being assumed as feasible.  The generalization to assuming known a strictly feasible matrix, but not necessarily the identity, is presented in section~\ref{sect.scaling}.  \vspace{1mm}

\begin{thm}  \label{thm.db}
Assume I is feasible for $ \mathrm{SDP} $.  $ \mathrm{NonSmoothed \, \,  Scheme} $  outputs $ Z $ which is feasible for $ \mathrm{SDP} $   and satisfies
\[ 
    \frac{\lin C, Z \rin - \mathrm{opt\_val}}{\lin C, I \rin - \mathrm{opt\_val}} \, \leq \, \epsilon \; . 
\] 
The total number of iterations of $ \mathrm{Subgradient \, \, Method} $  is bounded above by
\[ 
    \left( 9 \, \,  \mathrm{diam}^2  \, + 1 \right) \,   \cdot  \, \left( \,   \frac{1}{\epsilon^2}  \, + \,     \log_{3/2} \left(  \frac{\lin C, I \rin - \mathrm{opt\_val}} {\lin C, I \rin - \mathrm{val}_0 }  \right)    \right)  \; , 
\] 
where $ \mathrm{val}_0 :=  \lin C, U_0 \rin $ and $ U_0 $ is the input matrix. 
\end{thm} 
\noindent {\bf Proof:}  Proposition~\ref{prop.da}  shows the output matrix $ Y $ from Step 1 is feasible for SDP and satisfies
\[    \frac{\lin C, I \rin - \mathrm{opt\_val}}{\lin C, I \rin - \mathrm{val}} \, \geq \, \frac{1}{3} \; , \]
where $ \mathrm{val} = \lin C, Y \rin $.  Thus, by Corollary~\ref{cor.cb},  when $ X_0 = Y $ is input into Subgradient Method, along with $ R = \mathrm{diam} $ and $ N = \lceil (3 \, \mathrm{diam}/ \epsilon)^2   \rceil $, the projection $ Z(X)$ of the output $ X $ satisfies   
\[   \frac{\lin C, Z(X) \rin - \mathrm{opt\_val}}{\lin C, I \rin - \mathrm{opt\_val}} \, \leq \, \epsilon \; , \]
establishing correctness of NonSmoothed Scheme. 
\vspace{1mm}

  The bound for total iterations of Subgradient Method is immediate from the outer iteration bound of Proposition~\ref{prop.da},  and the choices for the number of iterations in Step 2 of NonSmoothed SubScheme and in Step 2 of NonSmoothed Scheme. \vspace{3mm} \hfill $ \Box $

  The following corollary is useful when an optimization problem is modeled so as to make $ I $ be on (or near) the central path. The proof follows standard lines in interior-point method theory, but nonetheless we include the proof for   \vspace{1mm} completeness.

\begin{cor}  \label{cor.dc}
 If $ I $ is on the central path and the input matrix is chosen as $ U_0 = I -  \smfrac{1}{ \lambda_{\max}(\pi(C))} \, \pi(C) $, then the same conclusions as in Theorem~\ref{thm.db}  apply but now with the number of $ \mathrm{Subgradient \, \, Method} $  iterations bounded above by
\[ 
     \left(  9  \, \,   \mathrm{diam}^2 \,  + 1 \right)  \,   \cdot \,  \left( \,   \frac{1}{\epsilon^2}  \, + \,     \log_{3/2}(n) \,   \right)     \; . 
\] 
\end{cor}
\noindent {\bf Proof:}  Assume $ I $ is on the central path, that is, assume for some $ \mu > 0 $ that $ I $ is the optimal solution for
\[ \begin{array}{rl}
   \min & \lin C, X \rin - \mu \, \ln ( \det(X)) \\
  \textrm{s.t.} & {\mathcal A}(X) = b \; . 
  \end{array} \]
Since the gradient of the objective function at $ X \succ 0 $ is $ C - X^{-1} $, a first-order optimality condition satisfied by $ I $ is that there exists a vector $ y $ for which
\[       C - I = {\mathcal A}^* y \; , \]
where $ {\mathcal A}^* y = \sum_i y_i A_i $ is the adjoint of $ {\mathcal A} $.  This implies that the projection of $ C $ and $ \mu I $ onto the nullspace of $ {\mathcal A} $ are identical.  Consequently, 
\begin{gather*}     \lin C, X \rin \leq  \lin C, I \rin \textrm{ and }  A(X) = b \\ \Leftrightarrow \\  \tr(X)  \leq  n \textrm{ and }  A(X) = b . 
\end{gather*} 
Hence, all $ X $ which are both feasible for SDP and have better objective value lie within the set $ \{ X \geq 0: \tr(X) \leq n \} $, a set which is contained within the ball of radius $ n $ centered at $ I $. Thus, all feasible $ X $ for SDP satisfy
\begin{align*} 
    \lin C, I \rin - \lin C, X \rin & = \lin \pi(C), I - X \rin \\
         & \leq \| \pi(C) \| \, n \; , 
\end{align*} 
that is,
\begin{equation} \label{eqn.dd}
         \lin C, I \rin - \mathrm{opt\_val} \leq n \| \pi(C) \| \; .
         \end{equation}

On the other hand, the feasible matrix $ U_0 = I -  \frac{1}{\lambda_{\max}(\pi(C))} \, \pi(C) $ lies distance at least $ 1 $ from $ I $, because the unit ball centered at $ I $ is contained in $ \Symp $ and because $ U_0 $ lies in the boundary of $ \Symp $. Hence,
\[   I - U_0 = \smfrac{\alpha }{\| \pi(C) \|}  \pi(C) \textrm{ for some $ \alpha \geq 1 $.}
\]
Consequently,
\begin{align}
  \lin C, I \rin - \lin C, U_0 \rin & = \lin \pi(C), I - U_0 \rin  \nonumber \\
                                    & = \alpha \| \pi(C) \| \nonumber \\
                                    & \geq  \| \pi(C) \| \; . \label{eqn.de}  
\end{align}

Combining (\ref{eqn.dd})  and (\ref{eqn.de}) gives
\[  \frac{\lin C, I \rin - \mathrm{opt\_val}}{\lin C, I \rin - \lin C, U_0 \rin } \, \leq \, n \; . \]
Substitution into Theorem~\ref{thm.db}  completes the proof. \vspace{3mm} \hfill $ \Box $

It is interesting to observe that for any fixed value of $ \epsilon $, if one is able to model a family of optimization problems as semidefinite programs $ \mathrm{SDP}(n) $ parameterized by $ n $ (the number of variables) in such a way that for every $ n $, both $ I_n $ is on (or near) the central path and for some $ p < 1/4 $, $ \mathrm{diam}(n) = O(n^p) $, then the iteration bound     
provided by the corollary is better than the best iteration bound established for interior-point methods, i.e., $ O( \sqrt{n}) $ iterations when $ \epsilon $ is fixed. As each iteration of Subgradient Method is cheap relative to the cost of an interior-point method iteration, in this case NonSmoothed Scheme wins hands down.  \vspace{1mm}

On the other hand, of course, if $ n $ is held fixed and $ \epsilon $ goes to zero, the bound  $ O(\log(1/\epsilon)) $ on the number of iterations for interior-point methods is massively better than the bound $ O(1/ \epsilon^2) $ for NonSmoothed Scheme.

\section{{\bf  Starting Points} $ \mathbf{E \neq I} $} \label{sect.scaling}

In this section it is observed that the theory and algorithms from previous sections are readily converted to the case that the starting point is a strictly-feasible matrix $ E \neq I $.  \vspace{1mm}

 As in the remarks closing section~\ref{sect.thy}, the relevant  concave function is
\begin{equation}  \label{eqn.ea}
      \lambda_{\min,E}(X) := \inf \{ \lambda \in \mathbb{R}: X - \lambda E \notin \Symp \} \; , 
      \end{equation} 
(the smallest eigenvalue of the matrix $ E^{-1/2} X E^{-1/2} $, where $ E^{1/2} $ is the positive definite matrix satisfying $ E = E^{1/2} E^{1/2} $).  As those remarks noted, the theory of that section is easily generalized, which for the present situation means replacing all occurrences of $ \lambda_{ \min} $ appearing in section~\ref{sect.thy}  by $ \lambda_{\min,E} $, while simultaneously replacing $ I $ by $ E $, assuming $ E $ is strictly feasible. \vspace{1mm}

The algorithms and theory from sections~\ref{sect.subgrad}  and \ref{sect.nonsmoothed}, however, are not so obviously extended.  At issue is that unlike $ X \mapsto \lambda_{\min}(X) $, the function $ X \mapsto \lambda_{\min, E}(X) $ need not be Lipschitz continuous with constant 1, a fact that was critical in the proof of Theorem~\ref{thm.ca}.  Thankfully, the issue is easily handled by changing from the trace inner product to the  inner product on $ \Symp $ used in interior-point method theory:
\[   \lin U, V \rin_E  \, :=  \,  \tr( E^{-1} U E^{-1} V) \; = \, \tr \big( (E^{-1/2}U E^{-1/2} ) \, (E^{-1/2}U E^{-1/2} ) \big)  \; .  \]

Thus, for example, when SubGradient Method computes a subgradient for iterate $ X_k $, the subgradient should be with respect to $ \langle \; , \; \rangle_E $ rather than with respect to $ \langle \; , \; \rangle $.  \vspace{1mm}

Likewise, when SubGradient Method projects a subgradient onto the subspace (\ref{eqn.cc}), the projection should be orthogonal with respect to $ \langle \; , \; \rangle_E $. \vspace{1mm}

Finally, the value $ \mathrm{diam} $ should be replaced by a value $ \mathrm{diam}_E $ satisfying   
\[  \mathrm{diam}_E \geq \max \{ \|X - Y \|_E: X, Y \in \mathrm{Level}_{\mathrm{val}} \} \quad  \textrm{ for all $ \mathrm{val} < \lin C, E \rin $}  \; .  \]

With these changes, all of the results of previous sections are valid with $ E $ in place of $ I $. \vspace{1mm}

For linear programming, the resulting changes to Subgradient Method are quickly described.  Letting $ e $ denote the known strictly-feasible point (not necessarily the vector of all ones), and letting $ \mathrm{val} $ be any value satisfying $ \mathrm{val} < c^T e $, the problem equivalent to LP is
\[ \begin{array}{rl}
   \max_x & \min_j x_j/e_j   \\
    \textrm{s.t.} & Ax = b \\
    & c^T x = \mathrm{val} \; . \end{array} \]
In applying Subgradient Method, the relevant inner product is
\[   \lin u, v \rin_e = \sum_j \frac{u_j v_j}{e_j^2} \; . \]
With respect to this inner product, the subgradients at $ x \in \mathbb{R}^n $ are the convex combinations of the vectors $ e'(k) $ for which $ x_k/ e_k = \min_j x_j/e_j $, where $ e'(k) $ has all coordinates equal to zero except for the $ k^{th} $ coordinate, which is equal to $ e_k $.  \vspace{1mm}

The $ \langle \; , \; \rangle_e $-orthogonal projections of the vectors $ e'(j) $ ($ j = 1, \ldots, n $) onto the nullspace of
\[   \bar{A} = \left[ \begin{matrix} A  \\ c^T \end{matrix} \right] \]
are the columns of the matrix $ \bar{P}_e $ that $ \langle \; , \; \rangle_e $-orthogonally projects $ \mathbb{R}^n $ onto the nullspace, that is,
\[   \bar{P}_e =  I - \Delta (e)^2 \, \bar{A}^T ( \bar{A}  \, \Delta(e)^2 \bar{A}^T)^{-1} \bar{A}  \; , \]
where $ \Delta (e) $ is the diagonal matrix with $ j^{th} $ diagonal entry equal to $ e_j $.  Thus, the projected subgradients relied upon by Subgradient Method are now the convex combination of the columns of $ \bar{P}_e $. \vspace{2mm}

Perhaps the easiest way to understand why all of the results of previous sections remain valid when $ I $ is replaced by $ E $ -- and $ \langle \; , \; \rangle $ is replaced by $ \langle \; , \; \rangle_E $ -- is to use the standard trick in the interior-point method literature of ``scaling'' SDP to an equivalent semidefinite program for which $ I $ is feasible. The equivalent semidefinite program (in variable $ Y $) is
\[ \left.  \begin{array}{rl}
 \min_Y & \lin E^{1/2} C E^{1/2}, Y \rin \\
 & {\mathcal A}(E^{1/2} Y E^{1/2}) = b \\
 & Y \succeq 0 \end{array} \right\} \, \mathrm{scaled\_SDP}
 \] 
The equivalence is seen by noting $ X $ is feasible for SDP if and only if $ Y = E^{-1/2} X E^{-1/2} $ is feasible for $ \mathrm{scaled\_SDP} $, and the objective value of $ X $ for $ \mathrm{SDP} $ is the same as the objective value of $ Y $ for $ \mathrm{scaled\_SDP} $. \vspace{1mm}

Moreover, the inner product $ \langle \; , \; \rangle_E $ is transformed into the trace inner product, in that $ \lin X_1, X_2 \rin_E = \lin Y_1, Y_2 \rin $ for $ Y_j = E^{-1/2} X_j E^{-1/2} $ ($ j = 1,2 $).
Thus, for example, the $ \langle \; , \; \rangle_E $-diameter of level set $ \mathrm{Level}_{\mathrm{val}}(\mathrm{SDP}) $ is the same as the $ \langle \; , \; \rangle $-diameter of level set $ \mathrm{Level}_{\mathrm{val}}(\mathrm{scaled\_SDP}) $. \vspace{1mm}

Lastly, as is straightforward but tedious to verify, each of the algorithms transforms as well.  For example, consider Subgradient Method, and fix two of the inputs, $ R $ and $ N $.  Assume the algorithm is applied with $ \langle \; , \; \rangle_E $ to solving the linearly-constrained problem equivalent to SDP:
\[ \begin{array}{rl}
    \max & \lambda_{\min,E}(X) \\
   \textrm{s.t.} & {\mathcal A}(X) = b \\
           & \lin C, X \rin = \mathrm{val} \; .  \end{array}  \]
Then $ X_0, X_1, \ldots, X_N $ is a possible resulting sequence if and only if $ Y_0, Y_1, \ldots, Y_N $ -- where $ Y_j := E^{-1/2} X_j E^{-1/2} $ -- is a possible resulting sequence when Subgradient Method is applied using the trace inner product to the problem equivalent to $ \mathrm{scaled\_SDP} $:  \vspace{1mm}
\[ \begin{array}{rl}
       \max & \lambda_{\min}(Y) \\
          \textrm{s.t.} & {\mathcal A}( E^{1/2} Y E^{1/2} ) = b \\
             & \lin E^{1/2} C E^{1/2}, Y \rin = \mathrm{val} \; . \end{array} \]

Applying any of the algorithms with strictly-feasible $ E $ -- and with inner product $ \langle \; , \; \rangle_E $ -- is, in other words, equivalent to scaling SDP, applying the algorithm with $ I $  and the trace inner product, and then unscaling the answer. \vspace{1mm}

We choose to assume $ I $ is feasible only to reduce notational clutter and make evident the simplicity of the main ideas. \vspace{1mm}

(Unfortunately, the trick of scaling does not generalize to hyperbolic programming, making the proofs in the forthcoming paper \cite{renegar2014hyperbolic}  necessarily more abstract than the ones here.) \vspace{1mm}

In closing the section, we remark that the results throughout the paper can be developed just as readily for semidefinite programs of, say, the form
\[ \begin{array}{rl}
       \min_x & c^T x \\
         \textrm{s.t.} & \sum_{i=1}^m x_i A_i  \succeq B \; . \end{array} \]
 One assumes known a strictly feasible point $ e' $, and relies upon the concave function $ X \mapsto \lambda_{\min,E}(X) $ specified in (\ref{eqn.ea}), letting $ E := \sum_{i=1}^m e_i' A_i \, - \, B $. The equivalent problem solvable by a subgradient method is
 \[ \begin{array}{rl}
    \max & \lambda_{\min, E}(\sum_{i=1}^m x_i A_i \, - \, B) \\
\textrm{s.t.} & c^T x = \mathrm{val} \; , \end{array} \]
for any value satisfying $ \mathrm{val} < c^T e' $.
The relevant inner product at $ e' $ is
 \[   \lin u, v \rin_{e'} :=  \lin \sum_{i=1}^m u_i A_i \, , \, \sum_{i=1}^m v_i A_i \rin_E \; . \]
 \vspace{1mm}

For the special case of a linear program
\[  \begin{array}{rl}
            \min & c^T x \\
   \textrm{s.t.} & A x \geq b \; , \end{array} \]
letting $ \alpha_i^T $ denote the $ i^{th} $ row of $ A $ and letting $ e := Ae' - b $, the equivalent problem for $ \mathrm{val} < c^T e' $ is  
\[ \begin{array}{rl}
    \max & \min_i \smfrac{1}{e_i} (\alpha_i^T x - b_i) \\
\textrm{s.t.} &  c^T x  = \mathrm{val} \; . \end{array} \]
 The relevant inner product is     
\[     \lin u, v \rin_{e'} = u^T A^T \,  \Delta(e)^{-2}\,  A v \; , \]
where $ \Delta(e) $ is the diagonal matrix with $ i^{th} $ diagonal entry equal to $ e_i $. \vspace{1mm}

\section{{\bf Corollaries for Nesterov's ``First First-Order Method''}}  \label{sect.nesterov}

Assume SDP has an optimal solution and $ I $ is feasible. (In exactly the same manner as previous results generalize to an arbitrary strictly-feasible initial matrix $ E $, so do all of the remaining results.)
\vspace{1mm}

Recall that given $ \epsilon > 0 $ and a value $ \mathrm{val} $ satisfying $ \mathrm{val} < \lin C, I \rin $, we wish to approximately solve the SDP equivalent problem
\begin{equation} \label{eqn.fa}
 \begin{array}{rl}
 \max & \lambda_{ \min}(X) \\
\textrm{s.t.} & {\mathcal A}(X) = b \\
      & \lin C, X \rin = \mathrm{val} \; , \end{array} 
      \end{equation}
where by ``approximately solve'' we mean that feasible $ X $ is computed for which 
\begin{equation}  \label{eqn.fb}
\lambda_{\min}(X^*_{\mathrm{val}}) - \lambda_{ \min}(X) \leq \epsilon' \;  
\end{equation} 
with $ \epsilon' $ satisfying
\[ 
 \epsilon' \leq \smfrac{\epsilon}{1 - \epsilon } \, \smfrac{\lin C, I \rin - \mathrm{val}   }{\lin C, I \rin - \mathrm{opt\_val} }  \; . 
\] 
Indeed, according to Proposition~\ref{prop.bd}, the projection $ Z = Z(X) $ will then satisfy
\[ 
\frac{\lin C, Z \rin - \mathrm{opt\_val}}{\lin C,I \rin - \mathrm{opt\_val}} \leq \epsilon  \; . 
\]
 \vspace{1mm}

With \cite{nesterov2005smooth}, Nesterov initiated a huge wave of research, by displaying that some significant non-smooth optimization problems can be efficiently solved by ``smoothing'' the problem and then applying optimal first-order methods for smooth functions.  In \cite{nesterov2007smoothing}, he extended the approach to include some problems within the domain of semidefinite programming.  Here he gave emphasis to the nonsmooth convex objective function $ X \mapsto \lambda_{ \max }(X) $, but the results trivially adapt to the concave function of interest to us, $ X \mapsto \lambda_{\min}(X) $.  \vspace{1mm}

For the nonsmooth concave function $ X \mapsto \lambda_{ \min}(X) $, the useful smoothing is
\[   f_{ \mu}(X) := - \mu \ln \sum_j e^{ - \lambda_j(X)/ \mu} \; , \]
where $ \mu > 0 $ is user-chosen, and where $ \lambda_1(X), \ldots, \lambda_n(X) $ are the eigenvalues of $ X $. For motivation, observe that for all $ X \in \Sym $,
\begin{equation}  \label{eqn.fe}
  \lambda_{\min}(X) - \mu \, \ln n \, \leq \,  f_{ \mu}(X) \, \leq \, \lambda_{\min}(X) \; . 
  \end{equation}
Not obvious, but which Nesterov proved,
\[ 
   \| \grad f_{ \mu }(X) - \grad f_{ \mu }(Y) \| \leq \smfrac{1}{\mu }  \| X - Y \| \; , 
\] 
where $ \| \, \, \| $ is the Frobenius norm.  That is, the gradient of $ f_{ \mu} $ is Lipschitz continuous, with constant $ 1/ \mu $. \vspace{1mm}
\vspace{1mm}

The smoothed version of (\ref{eqn.fa})   is
\begin{equation}  \label{eqn.ff}
 \begin{array}{rl}
 \max & f_{ \mu}(X) \\
  \textrm{s.t.} & {\mathcal A}(X) = b \\
  & \lin C, X \rin = \mathrm{val} \; . \end{array} 
  \end{equation}
Let $ X^*_{ \mathrm{val}}( \mu ) $ denote an optimal solution (which the reader should be careful to distinguish from $ X^*_{ \mathrm{val}} \; , $ an optimal solution for (\ref{eqn.fa})). \vspace{1mm}

In passing, we note that the gradient of $ f_{ \mu} $ at $ X $ is the
the matrix
\[   \grad f_{\mu }(X) = \smfrac{1}{\sum_j e^{ - \lambda_j(X)/ \mu}} \, \, Q \left[ \begin{smallmatrix} e^{-\lambda_1(X)/\mu } & & \\ & \ddots & \\ & & e^{-\lambda_n(X)/\mu } \end{smallmatrix} \right] Q^T , \]
where  $ X = Q \left[ \begin{smallmatrix} \lambda_1(X) & & \\ & \ddots & \\ & & \lambda_n(X)  \end{smallmatrix} \right] Q^T $ is an eigendecomposition of $ X $. \vspace{1mm}

  For the special case of linear programming, the function $ f_{\mu} $ becomes
  \[  f_{\mu}(x) = - \mu \ln \sum_j e^{ - x_j/ \mu} \; , \]
for which the gradient at $ x $ is the vector with $ j^{th} $ coordinate equal to
\[   \grad f_{\mu}(x)_j =  \frac{e^{ - x_j/ \mu}}{\sum_k e^{ - x_k/ \mu}} \; . \]

It is readily seen from (\ref{eqn.fe})  that for any value $ \epsilon' > 0 $ and for all $ X \in \Sym $, if $ \mu = \smfrac{1}{2} \epsilon'/\ln(n) $, then
\[  
    f_{ \mu} (X^*_{ \mathrm{val}}( \mu )) - f_{ \mu }(X) \leq \smfrac{1}{2} \epsilon' \quad \Rightarrow \quad \lambda_{\min}(X^*_{\mathrm{val}}) - \lambda_{ \min}(X) \leq \epsilon' \; . 
\] 
Consequently, in order to compute $ X $ which is feasible for (\ref{eqn.fa})  and satisfies (\ref{eqn.fb}), it suffices to fix  $ \mu = \smfrac{1}{2} \epsilon'/\ln(n) $ and compute $ X $ which is feasible for (\ref{eqn.ff})  and has objective value within $ \epsilon'/2 $ of $ f_{ \mu}( X^*_{\mathrm{val} }( \mu )) $.  \vspace{1mm}

Since the objective function in (\ref{eqn.ff}) is smooth and the only constraints are linear equations, we can apply many first-order methods.  It is only fitting that we rely on the original ``optimal'' first-order method for smooth functions, due to Nesterov \cite{nesterov1983method}, and which we refer to as ``Nesterov's first first-order method,'' or for brevity, ``Nesterov's first method.''\vspace{1mm}

  Letting $ X_0 $ denote an initial feasible point, then according to Theorem 2.2.2 in \cite{nesterov2004introductory}, Nesterov's first method produces a sequence of iterates $ X_0, X_1, \ldots $ satisfying
\[   f_{ \mu} ( X^*_{ \mathrm{val}}( \mu ) - f_{ \mu } (X_k) \leq \smfrac{4}{\mu \, (k+2)^2} \| X_0 - X^*_{ \mathrm{val}}( \mu )) \|^2  \]
(where we have used the fact that the Lipschitz constant for the gradient is $ 1/ \mu $).  Thus, $ f_{ \mu }(X_k) $ is within $ \epsilon'/2 $ of the optimal value if 
\[ k  \geq \sqrt{\frac{8}{ \mu \epsilon'}} \, \, \| X_0 - X^*_{ \mathrm{val}}( \mu )) \|   - 2  \; , \]
that is, if 
\[  
k \geq   \frac{4 \sqrt{\ln n} \,  \| X_0 - X^*_{ \mathrm{val}}( \mu )) \|}{ \epsilon'} - 2 \] 
(using  $ \mu = \smfrac{1}{2} \epsilon'/ \ln (n) $). We summarize these results in a theorem. \vspace{1mm}

\begin{thm} {\bf (Nesterov)} \label{thm.fa}
Let $ \epsilon' $ be any positive value, and $ \mu = \smfrac{1}{2} \epsilon'/ \ln(n) $. Assume $ X $ satisfies $ {\mathcal A}(X) = b $ and $ \lin C, X \rin = \mathrm{val} < \lin C, I \rin  $.  If Nesterov's first first-order method is applied to solving the smoothed problem (\ref{eqn.ff}), then the resulting iterates $ X_0 = X, X_1, X_2 , \ldots $ satisfy 
\[  
       k \geq  \frac{4 \sqrt{\ln n} \,  \| X - X^*_{ \mathrm{val}}( \mu )) \|}{ \epsilon'} - 2 \quad \Rightarrow \quad \lambda_{\min}(X^*_{\mathrm{val}}) - \lambda_{ \min} (X_k) \leq \epsilon' \; .\]       
\end{thm}
\vspace{3mm}

Recall that $ \mathrm{diam} $ is assumed to be a known quantity satisfying       
\[  \mathrm{diam} \geq \max \{ \|X - Y \|: X, Y \in \mathrm{Level}_{\mathrm{val}} \} \quad  \textrm{ for all $ \mathrm{val} < \lin C, I \rin $}  \; . \vspace{1mm} \]

\begin{cor} \label{cor.fb}
Assume $ X $ is strictly feasible for SDP and  $ \lin C, X \rin = \mathrm{val} < \lin C, I \rin $. 
 
For any value $ 0 < \epsilon' \leq 2 \, \lambda_{\min}(X) $, by letting $ \mu := \smfrac{1}{2} \epsilon'/\ln(n) $, if Nesterov's first first-order method is applied to solving the smoothed problem (\ref{eqn.ff}), then the resulting iterates $ X_0 = X, X_1, X_2,  \ldots $ satisfy
\[   k \geq \frac{4 \, \sqrt{\ln n} \, \,   \mathrm{diam}}{ \epsilon'} - 2 \quad \Rightarrow \quad \lambda_{\min}( X^*_{ \mathrm{val}}) - \lambda_{\min}(X_k) \leq \epsilon' \; . \]
\end{cor}
\noindent {\bf Proof:}  Assume $ X $, $ \epsilon' $ and $ \mu $ are as in the statement.  Then $ f_{\mu}(X) \geq 0 $, by (\ref{eqn.fe}) and $ \mu \leq \lambda_{\min}/ \ln(n) $.  Hence, $ f_{\mu}(X^*_{\mathrm{val}}(\mu)) \geq 0 $ (because $ f_{\mu}(X^*_{\mathrm{val}}(\mu)) \geq f_{\mu}(X) $). \vspace{1mm}

Also by (\ref{eqn.fe}),  $  
    Y \notin \Symp  \Rightarrow  f_{\mu}(Y) < 0 \; . $ Thus,  $ X^*_{\mathrm{val}}(\mu) $ is feasible for SDP. 
Hence, $   \| X - X^*_{\mathrm{val}}( \mu ) \| \leq \mathrm{diam} \; . $  
Substitution into Theorem~\ref{thm.fa}  completes the proof.   \vspace{3mm} \hfill $ \Box $

\begin{cor}   \label{cor.fc}
 Assume $ X $ is feasible for SDP and satisfies 
\begin{equation}  \label{eqn.fh}
  \lambda_{\min}(X) \geq \frac{1}{6}  \quad \textrm{and} \quad \frac{\lin C, I \rin - \mathrm{opt\_val}}{\lin C, I \rin - \mathrm{val}  } \leq 3 \; , 
  \end{equation} 
where $ \mathrm{val} := \lin C, X \rin $.   Let $ 0 < \epsilon < 1 $. 
  
For $ \mu = \smfrac{1}{6} \epsilon/ \ln(n) $, if Nesterov's first first-order method is applied to solving the smoothed problem (\ref{eqn.ff}), then the resulting iterates $ X_0 = X, X_1, X_2,  \ldots $ satisfy
\[  
    k \geq \frac{12 \, \sqrt{\ln n} \, \,   \mathrm{diam}}{ \epsilon} - 2 \quad \Rightarrow \quad \frac{\lin C, Z(X_k) \rin - \mathrm{opt\_val}}{\lin C, I \rin - \mathrm{opt\_val}} \, \leq \, \epsilon   \; .  
\] 
\end{cor}
\noindent {\bf Proof:}  
Let $ \epsilon' := \epsilon/3 \leq 1/3 $.  Then, by assumption, $ \epsilon' \leq 2 \lambda_{\min}(X) $, and hence Corollary~\ref{cor.fb} can be applied with 
\[  \mu = \smfrac{1}{2} \epsilon'/ \ln(n) = \smfrac{1}{6} \epsilon/ \ln(n) \; , \]
giving 
\begin{align*}
    k \geq \frac{12 \sqrt{\ln n} \, \, \mathrm{diam}}{ \epsilon} \, - \, 2 & \, \Rightarrow \, \lambda_{\min}(X^*_{\mathrm{val}}) - \lambda_{\min}(X) \, \leq \, \frac{\epsilon }{3} \\
   & \, \Rightarrow \,  
   \lambda_{\min}(X^*_{\mathrm{val}}) - \lambda_{\min}(X) \, \leq \frac{\epsilon }{1 - \epsilon} \, \, \frac{\lin C, I \rin - \mathrm{val}  }{\lin C, I \rin - \mathrm{opt\_val}} \; ,  
    \end{align*}
where the last implication is due to the rightmost inequality assumed in (\ref{eqn.fh}).  Invoking Proposition~\ref{prop.bd} completes the proof.   \hfill $ \Box $
\vspace{1mm}

\section{{\bf  The Smoothed Scheme}}  \label{sect.smoothed}

The presentation of the smoothed scheme is done in the same manner as the presentation of NonSmoothed Scheme in section~\ref{sect.nonsmoothed}, but now beginning with the following question: 
\begin{quote}
 Is it possible to efficiently move from an initial matrix $ U_0 $ satisfying $ {\mathcal A}(U_0) = b $ and $ \lin C, U_0 \rin < \lin C, I \rin $, to a matrix $ Y $ satisfying the conditions of Corollary~\ref{cor.fc}? 
 \end{quote}
Before providing an affirmative answer, for ease of reference we again display the pertinent optimization problems:
   \begin{equation}  \label{eqn.ga}
   \begin{array}{rl}
 \max & \lambda_{ \min}(X) \\
  \textrm{s.t.} &  {\mathcal A}(X) = b \\
   & \lin C, X \rin = \mathrm{val} \; , 
   \end{array} \end{equation} 
   \begin{equation}  \label{eqn.gb}
\begin{array}{rl}
 \max & f_{ \mu }(X) \\
\textrm{s.t.} & {\mathcal A}(X) = b \\
   & \lin C, X \rin = \mathrm{val} \; . 
   \end{array} \end{equation}
   Recall that $ X^*_{ \mathrm{val}} $ denotes any optimal solution for  (\ref{eqn.ga}) -- a problem equivalent to SDP (assuming $ \mathrm{val} < \lin C, I \rin $)  -- whereas $ X^*_{ \mathrm{val}}(\mu) $ denotes an optimal solution for (\ref{eqn.gb}) -- the smoothed version of (\ref{eqn.ga}).
   \vspace{1mm}

Consider the following computational procedure: 
\begin{itemize}

\item {\bf Smoothed SubScheme}
\begin{enumerate}

\addtocounter{enumi}{-1}

\item Initiation:

\begin{itemize}

\item Input:  A matrix $ U_0 $ that is feasible for SDP and satisfies both $ \lin C, U_0 \rin < \lin C, I \rin $ and $ \lambda_{\min}(U_0) = 1/6 $.  
\begin{itemize}

\item 
 If a matrix $ U $ is available satisfying only $ {\mathcal A}(U) = b $ and $ \lin C, U \rin < \lin C, I \rin $, then   $ U_0 = I + \smfrac{5}{6} \, \smfrac{1}{1 - \lambda_{\min}(U)} (U - I ) $ is acceptable input.
\end{itemize} 
\item Let $ \mathrm{val}_0 = \lin C, U_0 \rin $ 
\item Let $ \mu:= 1/(6 \, \ln n)  $ and $ N := \lceil 12 \, \sqrt{\ln n } \, \,  \mathrm{diam} - 2    \rceil $.
\item Let $ \ell = -1 $. 
\end{itemize}

\item Outer Iteration Counter Step:
\begin{itemize}

\item Let $ \ell + 1 \rightarrow \ell $

\end{itemize}
  
\item  Inner Iterations:

\begin{itemize}
\item Beginning at $ U_{ \ell} $, apply $ N $ iterations of Nesterov's first first-order method to the smoothed problem (\ref{eqn.gb}), with $ \mathrm{val} := \mathrm{val}_{\ell} $.  
\item  Let $ V_{ \ell} $ denote the final iterate.  

\end{itemize}

\item Check for Termination: 

\begin{itemize}

\item If $ \lambda_{\min}( V_{\ell}) \leq 1/3 $, then output $ Y = U_{ \ell} $ and terminate.  

\item Else, compute
\[  \qquad \qquad  \qquad  \qquad   U_{ \ell + 1} := I + \smfrac{5}{6} \, \smfrac{1}{1 - \lambda_{ \min}(V_{ \ell })} \, ( V_{ \ell} - I ) \; , \quad  \mathrm{val}_{\ell + 1} := \lin C, U_{ \ell + 1} \rin \; ,   \]
and go to Step 1.     
 
\end{itemize}

\end{enumerate}
\end{itemize}
\vspace{2mm}

\begin{prop}  \label{prop.ga}
$ \mathrm{Smoothed \, \, SubScheme} $   terminates with a matrix $ Y $ which is feasible for $ \mathrm{SDP} $   and satisfies
\begin{equation}  \label{eqn.gc}
  \lambda_{ \min }(Y) \,  = \, \frac{1}{6}  \; , \quad \frac{\lin C, I \rin - \mathrm{opt\_val}}{\lin C, I \rin - \mathrm{val}} \, \leq \, 3 \; , 
  \end{equation} 
where $ \mathrm{val} := \lin C, Y \rin $.  Moreover, the number of outer iterations does not exceed
\[  \log_{5/4} \left(  \frac{\lin C, I \rin - \mathrm{opt\_val}} {\lin C, I \rin - \mathrm{val}_0 }  \right) \;  \]
where $ \mathrm{val}_0 = \lin C, U_0 \rin $ and $ U_0 $ is the input matrix. 
\end{prop}
\noindent {\bf Proof:}  The proof -- especially the last half -- parallels that of Proposition~\ref{prop.da}.  Nonetheless, we include the proof in its entirety. \vspace{1mm}

It is easily verified that all of the matrices $ U_{ \ell} $ and $ V_{ \ell} $ computed by Smoothed SubScheme satisfy the SDP equations $ {\mathcal A}(X) = b $.  Moreover, $ U_{\ell} $ is feasible for SDP, because $ \lambda_{\min}(U_{ \ell}) = 1/6 $ (by construction). \vspace{1mm}

Fix $ \ell \geq 0 $ to be a value attained by the counter.  We now examine the effects of Steps 2 and 3. \vspace{1mm}

Applying Corollary~\ref{cor.fb}  with $ \epsilon' = 1/3 $ shows that in Step 2, the final iterate $ X_N $ computed by Nesterov's first method satisfies 
\[  \lambda_{\min}(X^*_{\mathrm{val}_{\ell}}) -  \lambda_{\min}(X_N) \leq 1/3 \; , \]
that is, 
\begin{equation}   \label{eqn.gd}
    \lambda_{\min}(X^*_{\mathrm{val}_{\ell}}) -  \lambda_{\min}(V_{\ell}) \leq 1/3 \; . 
    \end{equation}
    
Observe that
\begin{align*} 
\frac{\lin C, I \rin - \mathrm{val}_{ \ell} }{\lin C, I \rin - \mathrm{opt\_val}} & = 1 -  \lambda_{\min}(X^*_{\mathrm{val}_{\ell}}) \quad \textrm{(by Lemma~\ref{lem.bc})}  \\
  & \geq \smfrac{2}{3} - \lambda_{\min}( V_{ \ell}) \quad \textrm{(by (\ref{eqn.gd}))} \; . 
  \end{align*}
Hence, if the method terminates in Step 3 -- that is, if $ \lambda_{\min}(V_{\ell}) \leq 1/3 $ -- then the output matrix $ Y = U_{ \ell} $ satisfies
\[   \frac{\lin C, I \rin - \mathrm{val} }{\lin C, I \rin - \mathrm{opt\_val}} \, \geq \, \frac{1}{3} \; , \]
where $ \mathrm{val} :=  \lin C, Y \rin  = \mathrm{val}_{ \ell} $.  We have now verified that if the Smoothed SubScheme terminates, then the output $ Y $ is indeed feasible for SDP and satisfies the inequalities (\ref{eqn.gc}).  \vspace{1mm}

On the other hand, if the method does not terminate in Step 3, it computes the matrix $ U_{ \ell + 1} $ and its objective value, $ \mathrm{val}_{ \ell + 1} $.  Here, observe 
\begin{align*}
 \mathrm{val}_{ \ell + 1 }  & = \lin C, I \rin + \smfrac{5}{6} \, \smfrac{1}{1 - \lambda_{\min}(V_{\ell})} \, ( \mathrm{val}_{\ell} - \lin C, I \rin ) \\
  & \leq \lin C, I \rin - \smfrac{5}{4}  \, ( \lin C, I \rin - \mathrm{val}_{\ell} )  \; , 
  \end{align*}
because $ \mathrm{val}_{\ell} < \lin C, I \rin $ and  $ \lambda_{\min}(V_{ \ell}) \geq 1/3 $ (due to no termination). 
Hence,
\[   \frac{\lin C, I \rin - \mathrm{val}_{ \ell + 1} }{\lin C, I \rin - \mathrm{opt\_val}}  \geq \frac{5}{4} \, \, \frac{\lin C, I \rin - \mathrm{val}_{ \ell} }{\lin C, I \rin - \mathrm{opt\_val}} \; . \]
Since all values $ \mathrm{val}_{\ell} $ computed by the algorithm satisfy $ \mathrm{val}_{\ell} \geq \mathrm{opt\_val} $ (because $ U_{ \ell } $ is feasible for SDP),  it immediately follows that 
\[  \log_{5/4} \left(  \frac{\lin C, I \rin - \mathrm{opt\_val}} {\lin C, I \rin - \mathrm{val}_0 }  \right) \]
is an upper bound on the number of outer iterations.  \vspace{5mm} \hfill $ \Box $

Specifying our scheme based on Nesterov's first method, and analyzing  the scheme's complexity, both are now easily accomplished: \newpage
\begin{itemize}

\item {\bf Smoothed Scheme}

\begin{enumerate}

\addtocounter{enumi}{-1}

\item Input: A value $ 0 < \epsilon < 1 $, and $ U_0 \in \Sym $ satisfying $ {\mathcal A}(U_0) = b $ and $ \lambda_{\min}( U_0) = 1/6 $.   

\begin{itemize}

\item For example, the matrix $ U_0 = I - \smfrac{5}{6} \, \smfrac{1}{\lambda_{\max}(\pi(C))} \, \pi(C) $.

\end{itemize}
 
\item Apply Smoothed SubScheme with input $ U_0 $.  Let $ Y $ denote the output.

\item Beginning at $ Y $, apply $ \left\lceil \frac{12 \, \sqrt{\ln n} \, \, \mathrm{diam}}{ \epsilon }  - 2 \right\rceil$ iterations of Nesterov's first first-order method to the smoothed problem (\ref{eqn.gb}),  with $ \mathrm{val} :=  \lin C, Y \rin \; . $    Let $ X $ denote the output.

\item Compute and output the projection $ Z = Z(X) $, then terminate.   

\end{enumerate}

\end{itemize}

\begin{thm} \label{thm.gb}
Assume I is feasible.  $ \mathrm{Smoothed \, \, Scheme} $  outputs $ Z $ which is feasible for $ \mathrm{SDP} $   and satisfies
\begin{equation}  \label{eqn.ge}
    \frac{\lin C, Z \rin - \mathrm{opt\_val}}{\lin C, I \rin - \mathrm{opt\_val}} \, \leq \, \epsilon \; . 
    \end{equation} 
The total number of iterations of Nesterov's first first-order method is bounded above by
\[ 
   12 \, \sqrt{\ln n } \, \cdot \,   \mathrm{diam} \,  \cdot \,  \left(  \, \frac{1}{\epsilon} \, + \,       \log_{5/4} \left(  \frac{\lin C, I \rin - \mathrm{opt\_val}} {\lin C, I \rin - \mathrm{val}_0 }  \right) \, \right)  \; , 
\] 
where $ \mathrm{val}_0 :=  \lin C, U_0 \rin $ and $ U_0 $ is the input matrix. 
\end{thm} 
\noindent {\bf Proof:}  Proposition~\ref{prop.ga}  shows the matrix $ Y $ in Step 1 satisfies the conditions of Corollary~\ref{cor.fc}, which in turn shows the final output $ Z = Z(X) $ from Step 3 is feasible for SDP and satisfies (\ref{eqn.ge}). \vspace{1mm}

The bound on the total number of iterations of Nesterov's first  method is a straightforward consequence of the outer iteration bound from Proposition~\ref{prop.ga}, the choice for $ N $ in Smoothed SubScheme, and the number of iterates in Step 2 of Smoothed Scheme.   \vspace{3mm} \hfill $ \Box $

The proof of the following corollary proceeds exactly as does the proof of Corollary~\ref{cor.dc}. The added value 1 is due to $ \lambda_{\min}(U_0) = 1/6 $ in Smoothed Scheme -- unlike $ \lambda_{\min}(U_0) = 0 $  in NonSmoothed Scheme -- and 
\[  \log_{5/4}(\smfrac{1}{1 - \lambda_{\min}(U_0 )}) = \log_{5/4}(6/5) < 1 \; . \]  
 
\begin{cor}  \label{cor.gc}
If $ I $ is on the central path and $ U_0 = I - \smfrac{5}{6} \, \smfrac{1}{\lambda_{\max}(\pi(C))} \, \pi(C) $, then the same conclusions as in Theorem~\ref{thm.gb}  apply but now with the number of iterations of Nesterov's first first-order method bounded above by
\[  
12 \, \sqrt{\ln n } \, \cdot \,  \mathrm{diam}  \,  \cdot \,  \left(  \, \frac{1}{\epsilon} \, + \,      \log_{5/4}(n) \, + \,  1 \right)  \; . \]  
\end{cor}
\vspace{1mm}

\section{{\bf  Closing Remarks}} \label{sect.closing}

Similar to the observation immediately following Corollary~\ref{cor.dc}, we see from Corollary~\ref{cor.gc}  that for fixed $ \epsilon $, if one models a family of problems as semidefinite programs $ \mathrm{SDP}(n) $ where $ I_n $ is on (or near) the central path and for which there exists $ p < 1/2 $ satisfying $ \mathrm{diam}(n) = O(n^p) $, then Smoothed Scheme wins hands down over interior-point methods even on iteration count, let alone on total cost.     \vspace{1mm}

Interior-point methods, of course, win hands down if $ n $ is fixed and $ \epsilon $ goes to zero. \vspace{1mm}

For fixed $ n $, our iteration bound of $ O(1/ \epsilon) $ is much worse than the bound $ O(1/ \sqrt{\epsilon}) $ found in literature related to compressed sensing, where problems can be reduced to ones involving only the feasible regions $ \{ x \geq 0: \sum_j x_j = 1 \} $ -- or $ \{ X \succeq 0: \tr(X) = 1 \} $ -- for which tractable prox functions are known.  However, the approaches used there fail upon including  additional constraints, such as requiring $ X $ to satisfy a specific sparsity pattern, as is relevant in statistics for fitting a concentration matrix (the inverse of a covariance matrix) to empirical data, and as is relevant in some applications of semidefinite programming to combinatorial problems pertaining to graphs.  Among the obstacles is that tractable proximal operators remain unknown except for an extremely small universe of sets.   \vspace{1mm}

In this vein, we note that for the algorithms herein, imposing a specific sparsity pattern actually reduces work, assuming the known feasible matrix is $ E = I $.  Indeed, assuming the sparsity pattern is $ X_{ij} = 0 $ for $ (i,j) \in \mathrm{Zeros} $, and using $ \lin A_k, X \rin = b_i $ ($ k = 1, \ldots,  \ell ) $ to denote the remaining constraints, projecting a subgradient $ G $ is accomplished by overwriting by 0 the entries $ G_{i,j} $  for $ (i,j) \in \mathrm{Zeros}  $, and orthogonally projecting the resulting matrix onto the subspace
\[   \{ V: \lin \bar{C}, V \rin = 0 \textrm{ and }   \lin  \bar{A}_k, V \rin = 0 \textrm{ for $ k = 1, \ldots \ell  $} \}, \]
where $ \bar{C} $ (resp., $ \bar{A}_{k} $) is the matrix obtained by overwriting by 0 the $ (i,j)^{th} $ coordinate of $ C $ (resp., $ A_k $), for $ (i,j) \in \mathrm{Zeros} $.  If $ \ell $ -- the number of ``complicating constraints'' -- is small and the set $ \mathrm{Zeros} $ is large, this approach results in significant computational savings per iteration. \vspace{1mm}

Moreover, the resulting iteration cost is very cheap relative to the cost of an iteration of an interior-point method, where sparsity constraints must be handled like any other constraints, due to the inner product $ \langle \; , \; \rangle_{X_k} $ being dependent on the iterate, unlike first-order methods where the inner product $ \langle \; , \; \rangle = \langle \; , \; \rangle_I $ is held fixed throughout, an inner product for which sparsity constraints are ideally structured. \vspace{1mm}

Additional examples of the relevance of the algorithms and results, and their extensions to hyperbolic programming, are given in the forthcoming paper \cite{renegar2014hyperbolic}.

\bibliographystyle{plain}
\bibliography{first-order_algos_for_LP_and_SDP}

\end{document}